\documentclass[10pt]{amsart}

\usepackage[dvipsnames]{xcolor}
\usepackage{url,amssymb,graphicx,xypic,tikz}

\usepackage{hyperref}
\usepackage[nameinlink,capitalize]{cleveref}
\newcommand\myshade{85}
\colorlet{mylinkcolor}{violet}
\colorlet{mycitecolor}{YellowOrange}
\colorlet{myurlcolor}{Aquamarine}

\hypersetup{
  linkcolor  = mylinkcolor!\myshade!black,
  citecolor  = mycitecolor!\myshade!black,
  urlcolor   = myurlcolor!\myshade!black,
  colorlinks = true,
}

\newtheorem{theorem}{Theorem}[section]
\newtheorem*{theorem*}{Theorem}
\newtheorem{proposition}[theorem]{Proposition}
\newtheorem{corollary}[theorem]{Corollary}
\newtheorem{example}[theorem]{Example}
\newtheorem{lemma}[theorem]{Lemma}

\theoremstyle{definition}
\newtheorem{remark}[theorem]{Remark}

\newcommand{\A}{\mathcal{A}}
\newcommand{\CC}{\mathbf{C}}

\newcommand{\FI}{\mathtt{FI}}
\newcommand{\G}{\mathcal{G}}
\newcommand{\ZZ}{\mathbf{Z}}
\newcommand{\Res}{\mathrm{Res}}
\newcommand{\Ind}{\mathrm{Ind}}
\newcommand{\Lie}{\mathrm{Lie}}
\newcommand{\rk}{\mathrm{rk}}
\newcommand{\Sym}{\mathrm{Sym}}

\newcommand{\la}{\langle}
\newcommand{\ra}{\rangle}

\begin{document}
\begin{abstract}
  The group $G(m,1,n)$ consists of $n$-by-$n$ monomial matrices whose
  entries are $m$th roots of unity. It is generated by $n$ complex
  reflections acting on $\CC^n$. The reflecting hyperplanes give rise
  to a (hyperplane) arrangement $\G \subset \CC^n$. The internal
  zonotopal algebra of an arrangement is a finite dimensional algebra
  first studied by Holtz and Ron. Its dimension is the number of bases
  of the associated matroid with zero internal activity. In this paper
  we study the structure of the internal zonotopal algebra of the Gale
  dual of the reflection arrangement of $G(m,1,n)$, as a
  representation of this group. Our main result is a formula for the
  top degree component as an induced representation. We also provide
  results on representation stability, a connection to the Whitehouse
  representation in type~$\textup{A}$, and an analog of decreasing
  trees in type~$\textup{B}$.
\end{abstract}
\title[Internal Zonotopal Algebras and Reflection Groups]{Internal
  Zonotopal Algebras and the Monomial Reflection Groups $G(m,1,n)$}
\author{Andrew Berget}
\address{Western Washington Universty\\ Bellingham, WA USA}
\maketitle

\section{Introduction} A meta-problem in the theory of hyperplane
arrangements starts with an algebraic object $M(\A)$ derived from an
arrangement $\A$ and determines the extent to which the intersection
lattice $L(\A)$ determines $M(\A)$ up to isomorphism. The prototypical
example of this associates to an arrangement $\A \subset V$ the
cohomology ring of the complement $H^*(V - \A;\mathbf{Z})$. The
structure of this ring was determined by Orlik and Solomon \cite{OS1}
--- it is determined by the combinatorics of no-broken-circuit subsets
of the associated matroid.

There is an equivariant version of the meta-problem where the
arrangement $\A$ is fixed by the action of a group $W$ acting linearly
on $V$. When this happens one wants to determine the structure of
$M(\A)$ as a representation of $W$. The prototypical example here is
that of a real reflection group $W$ and its reflection arrangement
$\A \subset V$. One wants to understand the structure of
$H^*(V-\A;\mathbf{Z})$ as a representation of $W$. For complex
reflection groups this problem was first studied by Orlik and Solomon
\cite{OS2} and later by many others.

Of recent interest are the so-called zonotopal algebras of a
hyperplane arrangement
\cite{ardilaPostnikov1,ardilaPostnikov2,holtz-ron,lenzEJC,lenzIMRN}. The
zonotopal ideals of an arrangement $\A \subset V$ are ideals in
$Sym(V)$ generated by powers of elements of $v$. Specifically, the
$k$th zonotopal ideal of $\A$ is
\[
I_{\A,k}=\la h^{\max\{\rho_\A(h)+k+1,0\}} : h \in V \ra
\]
where $\rho_\A(h)$ is the number of hyperplanes in $\A$ that do not
contain $h$. The quotient $S_{\A,k} = Sym(V)/I_{\A,k}$ is the $k$th
zonotopal algebra of $\A$. Holtz and Ron \cite{holtz-ron} single out
the cases $k=-2,-1,0$ as being of particular interest, and call these
the internal, central and external zonotopal algebras of $\A$. The
internal case $k=-2$ exhibits dramatic subtleties in comparison to all
other cases $k >-1$ \cite{ardilaPostnikov2},
\cite[Section~7.4]{lenzEJC}. The initial motivation for studying the
internal zonotopal ideal was that its Macaulay inverse system consists
of those partial differential operators that leave continuous the box
spline associated to $\A$ (assuming $\A$ is unimodular)
\cite[Corollary~9]{lenzIMRN}.


In this paper we study the internal zonotopal algebra of certain
arrangements coming from complex reflection groups $W$. For
numerological reasons we do not study the naturally occurring
reflection arrangement of $W$, but instead its Gale dual. We are
specifically interested in the structure of this algebra as
representations of $W$.  We summarize our main results here.
\begin{theorem*}
  Let $\G$ be the reflection arrangement of one of the monomial groups
  $G(m,1,n) \subset GL_n(\CC)$, consisting of $n$-by-$n$ permutation
  matrices whose non-zero entries are $m$th roots of unity ($m>1$ and
  $n \geq 3$). Let $\G^\perp$ denote the Gale dual of $\G$. Then the
  degree $n-1$ component of $S_{\G^\perp,-2}$ is isomorphic to
  \[
    \Ind_C^W( \chi )
  \]
  as a representation of $W$. Here $C \subset W$ is the group
  generated by an $n$-cycle and the scalar matrix which multiplies by
  $e^{2\pi i/m}$, and $\chi$ is the character taking values
  $e^{2\pi i/n}$ and $e^{2\pi i (n-1)/m}$ at the respective generators
  of $C$.

  When $n$ and $m$ are coprime
  $\Ind_C^W( \chi ) \approx \Ind_{C'}^W( e^{2 \pi i/n} \cdot e^{2\pi
    i(n-1)/m})$, where $C'$ is the cyclic subgroup of $W$ generated by
  a Coxeter element.
\end{theorem*}
We complement this result with several others on the structure of the
internal zonotopal algebra $S_{\G^\perp,-2}$. We study the type A case
of the above theorem, when $m=1$. Based on our theorem above, it is
perhaps unsurprising that the Lie representation appears in type A,
however our perspective makes obvious the ``hidden'' action of a
larger symmetric group than is necessary on the Lie
representation, as observed by Mathieu \cite{mathieu} and Robinson and Whitehouse \cite{whitehouse}. For all $G(m,1,n)$, we give explicit generators of the
internal zonotopal ideal and then use this description to prove finite
generation in the sense of Sam and Snowden \cite{ss} and
representation stability as described by Gan and Li \cite{ganLi}. We
generalize a recurrence relation for the Whitehouse representation
which factorizes the regular representation of $G(m,1,n)$, $m>1$.  In
type B, when $m=2$, we use a Gr\"obner basis for the internal
zonotopal ideal to compute an analog of decreasing trees, which are
classic combinatorial objects.

When $m=2$ the main result bears some similarity to results of
N.~Bergeron \cite[Theorem~5.1]{bergeron} and, separately, Douglas
\cite{douglas}. Both authors were interested in the structure of the
cohomology ring of the complement of reflection arrangement of
$G(2,1,n)$. The character of this representation is also as induced
from the same subgroup described in our main theorem, albeit induced
from a different character. This \textit{a fortiori} similarity cannot
be explained by a general relationship between the internal zonotopal
algebra and the cohomology of the hyperplane complement.


\section{A worked example} Before we begin in earnest, we work a
complete example of the main theorem for the group $W =
G(2,1,3)$. This is the hyperoctohedral group of signed $3$-by-$3$
permutation matrices, which has $2^3 \cdot 3! = 48$ elements. $W$ acts
on $\CC^3$ in the natural way. The corresponding reflection arrangment
$\G \subset \CC^3$ is defined by the $9$ hyperplanes,
\[
  x_1 =  0,\quad
  x_2 = 0, \quad
  x_3 = 0,\quad
  x_1 \pm x_2 = 0,\quad
  x_1 \pm x_3 = 0,\quad
  x_2 \pm x_3 = 0.
\]
Here, $x_1,x_2,x_3$ are the standard basis vectors of $\CC^3$.

The Gale dual $\G^\perp$ of $\G$ is an arrangment of $9$ hyperplanes
in a $6$ dimensional vector space $K$, which has basis
\[
  y_{12}^0,y_{13}^0,y_{23}^0,
  y_{12}^1,y_{13}^1,y_{23}^1.
\]
The superscripts are read modulo $2$, the subscripts of the $y^0$'s
are anticommutative ($y_{ji}^0 = -y_{ij}^0$) and the subscripts of the
$y^1$'s are commutative ($y_{ij}^1 = y_{ji}^1$). The reflection group
$W$ acts on this vector space, where an honest permutation matrix
(i.e., ($0$/$1$)-matrix) permutes the subscripts of the $y$'s and the
matrix
\[
  \begin{bmatrix}
    -1 & 0 & 0 \\
    0  & 1 & 0 \\
    0  & 0 & 1
  \end{bmatrix}
\]
acts by fixing $y_{23}^0$ and $y_{23}^1$, while it sends $y_{12}^0$ to
$-y_{12}^1$ and $y_{13}^0$ to $-y_{13}^1$. This material is discussed
in \cref{sec:galeDuality} where it is shown that these rules really do
determine a representation of $W$.

The internal zonotopal ideal $I_{\G^\perp,-2}$ is quadratically generated in $Sym(K)$ by
\[
  (y_{12}^0)^2,
  (y_{12}^1)^2,\quad
  (y_{13}^0)^2,
  (y_{13}^1)^2,\quad
  (y_{23}^0)^2,
  (y_{23}^1)^2,\quad
  y_{12}^0y_{12}^1,\quad
  y_{13}^0y_{13}^1,\quad
  y_{23}^0y_{23}^1,\quad
\]
along with
\begin{align*}
  (y_{12}+y_{23} + y_{31})^2,  (-z_{12}+y_{23}+z_{31})^2 , (-z_{21} + y_{13}+ z_{32})^2 , (-z_{31}+y_{12}+z_{23})^2.
\end{align*}
This ideal is $W$-stable; see \cref{thm:ideal}. One can readily
compute with a computer algebra system like Macualay2 \cite{M2} (which
we used extensively to conjecture our results) that the Hilbert series
of the quotient $Sym(K)/I_{\G^\perp,-2}$ is
$1+6q+8q^2 = (1+2q)(1+4q)$. This follows from \cref{cor:hilbG}. A
Gr\"obner basis for this ideal is studied in \cref{sec:typeB}.

From the Hilbert series one readily obtains that the largest degree
component of $Sym(K)/I_{\G^\perp,-2}$ that is non-zero is an
$8$-dimensional representation of $W$. The main theorem of this paper,
\cref{thm:ind}, identifies this character as an induced character,
namely,
\[
  \Ind_C^W( \chi)
\]
where $C$ is the abelian subgroup of $W$ with generators
\[
  \begin{bmatrix}
    0 & 1 & 0 \\
    0 & 0 & 1 \\
    1 & 0 & 0
  \end{bmatrix},
  \quad
  \begin{bmatrix}
    -1 & 0 & 0 \\
    0 & -1 & 0 \\
    0 & 0 & -1
  \end{bmatrix}
\]
and $\chi$ takes the values $e^{2\pi i /3}$ and $(-1)^{3-1} = 1$,
respectively, at the generators of $C$.  Since $C$ is an abelian group
of order $6$ it must be that $C$ is cyclic, and the cyclic generator
must be Coxeter element $c \in W$. It follows that the $8$ dimensional top
of $Sym(K)/I_{\G^\perp,-2}$ is isomorphic to
$\Ind_{\langle c \rangle}^W( e^{2\pi i/3})$ as a representation of
$W$. This is \cref{cor:coxeter}.

The group $G(2,1,2)$ is a subgroup of $G(2,1,3)$ by viewing the former
as signed $3$-by-$3$ permutation matrices of the form
\[
  \begin{bmatrix}
    * & * & 0 \\
    * & * & 0 \\
    0 & 0 & 1
  \end{bmatrix}
\]
There are $8$ such matrices. We may thus view the $8$ dimensional top
of $Sym(K)/I_{\G^\perp,-2}$ as a representation of
$G(2,1,2)$. \cref{cor:cyclic} identifies this as the regular
representation of $G(2,1,2)$.


\section{Representation theory background} In this section we briefly
recall some notions from representation theory. All the material
needed in this paper can be found in, e.g., \cite[Chapters
1--3]{fultonHarris}.

A representation of a finite group $G$ is a group homomorphism
$G \to GL(V)$, where $V$ is a finite dimensional vector space (say,
over the complex numbers). The simplest example of a representation of
$G$ is the group algebra $\mathbf{C}G$, whose basis elements are the
elements of $G$ and where $G$ acts on the left by multiplication in
the group. The group algebra $\mathbf{C}G$ is a finite dimensional
algebra and a representation of $G$ is equivalent to a (left) module
over $\mathbf{C}G$. We conflate the two notions freely.

Let $G$ be a group and $H \subset G$ a subgroup. Given a
representation $U$  of $H$ one may
extend scalars to obtain a representation of $G$:
\[
  U\quad \leadsto\quad \Ind_H^G(U) := \mathbf{C}G
  \otimes_{\mathbf{C}H} U.
\]
This process is called induction, and $\Ind_H^G(U)$ is called an
induced representation. The dimension of $\Ind_H^G(U)$ is
$ \dim(U) \cdot |G|/|H| $. The associativity of tensor products proves
that there is a natural isomorphism of representation of $F$,
\begin{equation}\label{eq:ind}
\Ind_G^F(\Ind_H^G (U)) \approx \Ind_H^F(U).
\end{equation}
This is refered to as the transitivity of induction.

Say that $U$ is a representation of $H$ and is contained in a
representation $V$ of $G$. By definition, there is a map
$\mathbf{C}G \otimes_{\mathbf{C}H} U \to V$. If the $G$ orbit of $U$
is all of $V$ then the above map is surjective. If, in addtion, the
dimension of $V$ is exactly $\dim(U) \cdot |G|/|H| $, the above map is
an isomorphism and hence $V \approx \Ind_H^G(U)$. In the sequel we
will identify induced representation by demonstrating these
properties.

The character of a representation $V$ of $G$ is the function
$\chi_V : G \to \mathbf{C}$, $g \mapsto trace(g : V \to V)$. This
function determines $V$ up to isomorphism, and is constant on
conjugacy classes of $G$. When $G$ is
cyclic it suffices to describe a one dimensional representation by its
character value at a generator, as we will often do.  The character of
$\Ind_H^G(U)$ is easy to determine from the character of $U$. The
induced character at $g$ is zero unless $g$ is conjugate to an element
of $H$. On elements $h \in H$, the induced character is equal to
$\frac{|G|}{|H||h^G|}\sum_{x \in H \cap h^G} \chi_U(x)$. Here $h^G$ is
the conjugacy class of $h$.

\begin{example}[Stanley {\cite[Lemma~7.2]{stanley}}]\label{ex:lie}
  Let $C$ be the cyclic subgroup generated by an $n$-cycle $c$ in
  $\mathfrak{S}_n$. Abusing notation, let $e^{2\pi i/n}$ denote the one
  dimensional representation of $C$ where $c$ acts by $e^{2\pi i
    /n}$. Then, the character of $\Ind_C^{\mathfrak{S}_n} (e^{2\pi i /n})$
  at $c^d$, for $d$ a divisor of $n$, is
  \[
    \frac{n!}{n \cdot (n!/(n/d)^d d! )} \sum_{\substack{1 \leq k \leq
        n \\ k \textup{ prime to }n/d}} e^{2\pi i dk/n} = (d-1)!
    (n/d)^{d-1} \mu(n/d).
  \]
  Here $\mu$ is the usual number-theoretic M\"obius function.
\end{example}

\section{Zonotopal Algebras}\label{sec:zonotopalAlgebra}
Let $\A$ be a central arrangement of hyperplanes in a complex finite
dimensional vector space $V$. We let $M=M(\A)$ denote the matroid of
$\A$. Recall that this is the simplicial complex on $\A$ whose faces
(alias independent sets) are collections of hyperplanes whose defining
linear forms are linearly independent. Maximal faces in $M$ are
referred to as its bases.

Define a function $\rho_\A: V \to \mathbf{N}$ by the rule
\[
\rho_\A(h) = \textup{ the number of hyperplanes in $\A$ not containing }h,
\]
and use this to define the ideal
\[
I_{\A,k} = \la h^{\max\{\rho_{\A}(h)+k+1,0\}} : h \in V \ra \subset Sym(V).
\]

The quotient $S_{\A,k} := Sym(V)/I_{\A,k}$ has Krull dimension zero,
and so it is a finite dimensional complex vector space.  The ring
$Sym(V)$ is graded in the usual way, and we denote the $m$th graded
piece of a graded module over this ring by $(-)_m$. The ideal
$I_{\A,k}$ is graded, since it is generated by powers of linear forms,
which are homogeneous. The \textit{Hilbert series} of $S_{\A,k}$ is
the generating function for the dimensions of the graded pieces of
this graded $Sym(V)$-module. Since $(S_{\A,k})_m = 0$ for $m$
sufficiently large there is a biggest $m$ for which
$(S_{\A,k})_m \neq 0$ and we call this the \textit{top} of quotient
and denote it by $S_{\A,k}(\textrm{top})$ to ease the proliferation of
superscripts and subscripts.

\begin{example}\label{ex:zonotopalAlg1}
  Let $\A$ be the arrangement in $\mathbf{C}^2$ whose hyperplanes are
  defined by $x=0$, $y=0$ and $x \pm y = 0$. The real picture is shown
  below by the black lines.
\[
\begin{tikzpicture}[baseline=0pt]  
  \draw[very thick,->,purple] (0,0)--(1/2,.8);
  \node[above right, purple] at (1/2,.8) {$h_1$};
  \node[below right, blue] at (0,-.5) {$h_2$};
  \draw[thick] (-1,1)--(1,-1) ;
  \draw[thick] (-1,0)--(1,0);
  \draw[thick] (0,-1)--(0,1);
  \draw[thick] (-1,-1)--(1,1); 
  \draw[very thick,->,blue] (0,0)--(0,-.5);
\end{tikzpicture}
\]
Then, $h_1^3 \in I_{\A,-2}$ and $h_2^2 \in I_{\A,-2}$. One computes
that $\langle x^2,y^2, (x\pm y)^2 \rangle \subset I_{\A,-2}$.
\end{example}

\subsection{Hilbert series} We will describe the Hilbert series of
$S_{\A,k}$ in terms of the matroid of $\A$. To do so, we need the
Tutte polynomial of a matroid. We take the most expedient route. Given
a matroid $M$ with ground set $E$, the \textit{rank} of a subset $S$
of $E$, is size of a largest independent set of $M$ contained in
$S$. The \textit{Tutte polynomial} of a matroid $M$ with ground set
$E$ and rank function $\rk : 2^E \to \mathbf{N}$ is the bivariate
polynomial
\[
T_M(p,q) = \sum_{A \subset E} (p-1)^{\rk(E)-\rk(A)}(q-1)^{|A|-\rk(A)}.
\]
The Tutte polynomial of a matroid is universal in the sense that any
matroid invariant taking values in an abelian group satisfying a
generalized deletion-contraction identity must be an evaluation of the
Tutte polynomial.

For $k \geq -2$ there is a short exact sequence relating $S_{\A,k}$ to
$S_{\A \setminus H,k}$ and $S_{\A/H,k}$, where $\A \setminus H$ is
$\A$ with one of its defining hyperplanes $H$ removed, and $\A/H$ is
the arrangement obtained by intersecting $\A$ with $H$. When interpreted
at the level of Hilbert series this becomes a deletion-contraction
relation and we have the following result.
\begin{theorem}\label{thm:ap}
  Let $\A$ be a central, essential arrangement of hyperplanes in $V$.
  Let $M=M(\A)$ denote the matroid of the arrangement $\A$, which
  consists of $m$ hyperplanes.  The Hilbert series of $S_{\A,k}$ is
  equal to
  \begin{enumerate}
  \item $q^{m-\rk(M)} T_M(1+q,1/q)$ if $k=0$;\label{external}
  \item $q^{m-\rk(M)} T_M(1,1/q)$ if $k=-1$;\label{central}
  \item $q^{m-\rk(M)} T_M(0,1/q)$ if $k=-2$.\label{internal}
  \end{enumerate}
\end{theorem}
The algebras $S_{\A,k}$ occurring in the theorem are, respectively,
referred to as the external, central and internal zonotopal algebras
of $\A$ by Holtz and Ron.  Parts \eqref{external} and \eqref{central}
were discovered by many authors in many contexts (we mention only
\cite{ardilaPostnikov1,ardilaPostnikov2,holtz-ron}, but the interested
reader should look to these for more extensive references). The
internal zonotopal algebra is more subtle than its central and
external counterparts, and was not discovered until the work of Holtz
and Ron \cite{holtz-ron} where they proved \eqref{internal}.

In the external and central cases, bases of the Macaulay inverse
system of the associated ideals are given by certain products of
linear forms defining $\A$. Analogous results were conjectured to hold
in the internal case and an incorrect proof was given in
\cite{ardilaPostnikov1}. This was corrected in
\cite{ardilaPostnikov2}, and a combinatorial basis was later given by
Lenz \cite{lenzIMRN} when $\A$ is unimodular. For non-unimodular $\A$,
there is no known canonical basis of the Macaulay inverse system of
$I_{\A,-2}$ described by the matroid of $\A$. This is particularly
striking since the dimension of $S_{\A,-2}$ is $T_M(0,1)$, a very well
known number: it enumerates the number of bases of $M$ with internal
activity zero and is the reduced Euler characteristic of the
independent set complex of $M$. Finding a combinatorial basis for
$S_{\A,-2}$ for general $\A$ is a tantalizing open problem.

A smaller generating set for $I_{\A,k}$ suffices than the one given
when $k \in \{-2,-1,0\}$. To describe the smaller generating set, we
recall that a \textit{line} of the arrangement $\A$ is a
one-dimensional intersection of hyperplanes in $\A$.
\begin{theorem}[Ardila--Postnikov \cite{ardilaPostnikov1}]\label{thm:line}
  Let $\A$ be a central, essential arrangement of hyperplanes in
  $V$. The ideal
\[
I'_{\A,k} = \la h^{\rho_{\A}(h)+k+1} : h \textup{ spans a line of }\A \ra
\]
is equal to $I_{\A,k}$ for $k \in \{-2,-1,0\}$. 
\end{theorem}
\begin{example}
  Continuing \cref{ex:zonotopalAlg1}, the containment
  $\langle x^2,y^2, (x\pm y)^2 \rangle \subset I_{\mathcal{A},-2}$ is
  an equality by the previous result. The Tutte polynonmial of $\A$ is
  $p^2+q^2+2p+2q$, by definition. The Hilbert series of $S_{\A,-2}$ is
  $1+2q$, as can be found by the computation of the ideal, or using
  \cref{thm:ap}.
\end{example}

\section{Complex Reflection Groups}\label{sec:reflGps} A comprehensive treatment of
reflection groups in unitary spaces is Lehrer and Taylor's book
\cite{lehrerTaylor}.  Let $V$ be a complex vector space. A
(generalized) reflection is an element of $GL(V)$ that fixes a
hyperplane point-wise and has finite order.

A complex reflection group is a finite subgroup of $GL(V)$ generated
by reflections. Such groups have been classified by Shephard and Todd,
and Serre and Chevallay. There appears a single infinite family of
groups $G(me,e,n)$, $m,e,n \geq 1$, called the monomial groups as well
as $34$ exceptional groups. The monomial group $G(me,e,n)$ consists of
$n$-by-$n$ permutation matrices whose entries are $me$-th roots of
unit, the product of which is a $e$-th root of unity.

Associated to a complex reflection group $W$ is its reflection
arrangement $\mathcal{A} \subset V$. This is the arrangement of
hyperplanes that are fixed by all the reflections $W$.

The reflection arrangement $\mathcal{G} = \mathcal{G}_{m,1,n}$
associated to the group $G(m,1,n)$ ($m>1$) is defined in
$V = \mathbf{C}^n$ by the coordinate hyperplanes $x_i = 0$
($1 \leq i \leq n$) and the hyperplanes
\[
x_i - e^{2\pi i k/m} x_j \qquad (1 \leq i < j \leq n,\ 1 \leq k \leq
m).
\]

\begin{example}
  The arrangement in \cref{ex:zonotopalAlg1} is the reflection
  arrangement of $G(2,1,2)$.  The group $G(2,1,2)$ is generated by the
  two matrices displayed below.
  \[
    \begin{bmatrix}
      0 & 1\\ 1& 0 
    \end{bmatrix},
    \begin{bmatrix}
      -1 & 0 \\ 0 & 1
    \end{bmatrix} \in G(2,1,2).
  \]
\end{example}

The group $G(1,1,n)$ is the symmetric group, which we will usually
write as $\mathfrak{S}_n$. There is an obvious injection
$\mathfrak{S}_n \to G(m,1,n)$, a fact we will make often use of. The
group $G(m,1,n)$ can also be thought of as the wreath product
$\mathbf{Z}_m \wr \mathfrak{S}_n$.
\subsection{Invariants and degrees}
The Chevalley--Shephard--Todd theorem characterizes complex reflection
groups as those subgroups $W \subset GL(V)$ for which the ring of
polynomial invariants $Sym(V^*)^W$ is itself a polynomial
ring. Algebraically independent generators of the invariant ring will
not be unique, but their degrees $d_1 \leq d_2 \leq \dots \leq d_\ell$
will be and are referred to as the \textit{degrees of $W$}.

The Coxeter number of $W$ is the largest degree of $W$ and will be
denoted by $h$. For reflection groups $W \subset GL(V)$ generated by
$\dim(V)$ reflections (which are called well-generated) the integers
$d_i'$ satisfying $d_i + d'_{n-i+1}= h $ are called the \textit{codegrees
of $W$}. 

The following result can be found in \cite[Appendix~D.2]{lehrerTaylor}.
\begin{proposition}
  The degrees of $G(m,1,n)$ are $m,2m,3m,\dots,mn$. The codegrees of
  $G(m,1,n)$ are $0,m,2m,\dots,(n-1)m$.
\end{proposition}

The Tutte polynomial evaluation
$\chi_\A(q) = (-1)^{\rk(\A)}T_\mathcal{A}(1-q,0)$ is called the
characteristic polynomial of the arrangement. The polynomial
$(-q)^{\rk(\A)} \chi_\A(-1/q)$ is Hilbert series of the cohomology
ring of the complement of the arrangment $V -\A$. When $\A$ is a
reflection arrangement the relationship to the codegrees of $W$ is
this.
\begin{theorem}[Orlik--Solomon {\cite[Theorem~5.5]{OS2}}]\label{thm:os}
  Let $\A$ be the reflection arrangement of a well-generated complex
  reflection group $W$ with codegrees $d_1',\dots,d'_n$. Then,
  \[
    q^{\rk(\A)} T_\mathcal{A}(1+1/q,0) = \prod_{i=1}^n (1+(1+d'_i) q)
  \]
  is the Hilbert series of $H^*(V-\A;\mathbf{Z})$.
\end{theorem}

\begin{example}
  In our running example, the Tutte polynomial of $\A$ is
  $p^2+q^2 + 2p+2q$. We see that 
  \[q^{\rk(\A)}T_\A(1+1/q,0) = q^2 ((1+1/q)^2 + 2(1+1/q)) 
    = (1+q)(1+3q).
  \]
  The codegrees of $G(2,1,2)$ are $0$ and $2$. 
\end{example}


\section{Gale Duality} \label{sec:galeDuality} 
\subsection{Definition and properties}
Any arrangement of hyperplanes has a dual arrangement, determined as
follows. If $\mathcal{A} \subset V$ is an arrangement then we can
construct $\mathbf{C}^{\A}$, a vector space with a basis is in
bijection with the hyperplanes defining $\mathcal{A}$.

Any choice of linear functionals defining $\A$ determines a linear map
$\CC^\A \to V^*$. The kernel of this map is denoted $K$, so we have an
exact sequence
\[
0 \to K \to \CC^\A \to V^*.
\]
Dualizing gives an exact sequence
\[
0 \leftarrow K^* \leftarrow (\CC^\A)^* \leftarrow V.
\]
Since there is a natural basis of $(\CC^\A)^*$ indexed by the
hyperplanes in $\A$, the map $(\CC^\A)^* \to K^*$ determines an
arrangement in $K$. This is the \textit{Gale dual of $\A$}, which we
denote by $\A^\perp$. Assuming that the map to $V^*$ above is
surjective, Gale duality is an honest duality in the sense that
$(\A^\perp)^\perp = \A$.

Gale duality can be interpreted at the level of matroids. Let $M$ be a
matroid whose ground set we denote by $E$ and whose bases are denoted
$\mathcal{B}(M)$. The \textit{dual matroid $M^\perp$} of $M$ has ground set $E$
and
\[
\mathcal{B}(M^\perp) =
\{
E-B : B \in \mathcal{B}(M).
\}
\]
If $M$ is the matroid associated to an arrangement $\A$ then $M^\perp$
is the matroid associated to $\A^\perp$.
\begin{proposition}[{\cite[Proposition~6.2.4]{oxley}}]\label{prop:dual}
  Let $M$ be a matroid and $M^\perp$ its dual matroid. Then the Tutte polynomials of $M$ and $M^\perp$ are related by 
  \[
    T_{M^\perp}(x,y) = T_M(y,x).
  \]
\end{proposition}
\begin{example}
  In our running example, $\A$ is self-dual. In general, to find the
  coordinates of the Gale dual of an arrangement whose normal vectors
  are the columns of the block matrix $[I|A]$, where $I$ is an
  appropriately sized identity matrix, one computes $[-A^t|I]$. In the
  present example we have
  \[
    \begin{bmatrix}
      1 & 0 & 1 & 1 \\
      0 & 1 & 1 & -1
    \end{bmatrix}\leadsto 
    \begin{bmatrix}
      -1 & -1 & 1 & 0 \\
      -1 &  1 & 0 & 1 
    \end{bmatrix},
  \]
  and after relabeling and scaling the columns, we get the same
  arrangment. The Tutte polynomial of $\A$ is $p^2+q^2 + 2p + 2q$,
  which is visibly symmetric in $p$ and $q$.
\end{example}
We employ the previous result to compute Hilbert series as such.
\begin{proposition}
  Let $\A \subset V$ be the reflection arrangement of well-generated
  complex reflection group $W$. Say that the codegrees of $W$ are
  $d'_1,\dots,d'_n$ Then the Hilbert series of $S_{\A^\perp,-2}$ is
  \[
    \prod_{i=1}^n (1+d'_i q)
  \]
  In particular, the dimension of the degree $n-1$ piece of
  $S_{\A^\perp,-2}$ is the product of the non-zero codegrees of $W$.
\end{proposition}
\begin{proof}
  It follows from \cref{thm:ap} that the Hilbert series in
  question is the Tutte polynomial evaluation
  $q^nT_{\A^\perp}(0,1/q)$. By \cref{prop:dual} this is
  $q^n T_{\A}(1/q,0)$ and by the \cref{thm:os} we obtain
  \[
    q^n T_\A(1/q,0) = \prod_{i=1}^n(1+d'_i q).\qedhere
  \]
\end{proof}
\begin{corollary}\label{cor:hilbG}
  Let $\G$ be the reflection arrangement of the group $G(m,1,n)$ and
  let $\mathcal{G}^\perp$ be its dual. Then the Hilbert series of
  $S_{\G^\perp,-2}$ is
  \[
    (1+mq)(1+2m q)\dots (1+ (n-1)m q).
  \]
\end{corollary}
\begin{proof}
  This follows from the previous result since the codegrees of
  $G(m,1,n)$ are $0,m,2m,\dots,(n-1)m$.
\end{proof}

\subsection{Group actions and duality}
Assume that $V$ is a representation of a group $W$ and that $W$
stabilizes $\A$ in the sense that for any $w \in W$, $w \A = \A$. The
arrangement $\A^\perp$ can be made to carry an action of $W$ in a
natural way, as we now explain. Using the natural action of $W$ on
$V^*$, the elements of $W$ act by permuting and rescaling of the
linear functionals defining $\A$. This affords an action of $W$ on
$\CC^\A$, whereby $W$ acts by generalized permutation matrices. This
action is cooked up so that the map $\CC^\A \to V^*$ is equivariant.

Since $K \subset \CC^\A$ is the kernel of a map of $W$-modules it is
$W$ invariant. Since the map $(\CC^\A)^* \to K^*$ is a map of
$W$-modules we see that the arrangement $\A^\perp$ is fixed in $K$ by
the natural action of $W$.
\subsection{The Gale dual of the reflection arrangement of $G(m,1,n)$}\label{sec:galedualG}
Let $\mathcal{G} \subset \CC^n$ denote the reflection arrangement
associated to the group $G(m,1,n)$. We assume that $m > 1$ until
\cref{sec:typeA}. In this section we give an explicit
coordinatization of the Gale dual $\mathcal{G}^\perp$ and describe how
$W = G(m,1,n)$ acts on the dual space $K$ in which $\mathcal{G}^\perp$
lives.

We begin by labeling the hyperplanes that define $\G$. To ease the
notation (only slightly) we will write $\omega$ for $e^{2\pi i/m}$ Let
$h_i$ denote the hyperplane defined by the vanishing of the $i$-th
coordinate function $e_i = 0$ in $\CC^n$. Let $h_{ij}^k$ denote the
hyperplane defined by $e_i - \omega^k e_j = 0$, $i<j$. We will use the
convention that the superscript in this notation is read modulo $m$
and that
\[
h_{ji}^k = -\omega^k h_{ij}^{-k}.
\]
The $h$'s form the basis for the vector space $\CC^\G$ and the linear
map $\CC^\G \to \CC^n$ sends $h_i \mapsto e_i$ and
$h_{ij}^k \mapsto e_i -\omega^k e_j$, where we are identifying $\CC^n$
with its dual space.

The group $W$ is generated by the $n$-by-$n$ permutation matrices,
which we identify with the symmetric group $\mathfrak S_n$, and the diagonal
matrices $g_i$ that have all diagonal entries $1$ except the $i$-th
which is $\omega$.  The action of $W$ on $\CC^\G$ is described in the
following proposition.
\begin{proposition}
  Given $\pi \in \mathfrak S_n \subset W$, 
  \[
    \pi h_i = h_{\pi (i)}, \qquad \pi h_{ij}^k  = h_{\pi(i)\pi(j)}^k.
  \]
  For all $i,j$,
  \[
    g_i h_i = \omega h_i, \qquad
    g_i h_{ji}^k = \omega h_{ij}^{k-1}, \qquad g_j h_{ij}^k  = h_{ij}^{k+1}.
  \]
  If $\ell \notin \{i,j\}$ then $g_\ell h_{ij}^k = h_{ij}^k$.
\end{proposition}

The kernel $K$ of the map $\CC^\G \to \CC^n$ is $m \binom{n}{2}$
dimensional and has basis given by
\[
  y_{ij}^k = (h_i - \omega^k h_j) - h_{ij}^k, \quad (1 \leq i < j \leq n, 0 \leq k < m).
\]
As before, we define $y_{ji}^k = -\omega y_{ij}^{-k}$ and read the
supersscript $k$ modulo $m$.  The arrangement $\G^\perp \subset K$
comes from restricting the linear functionals dual to the $h$'s to
$K$. The action of $W$ on $K$ is immediate from the previous
proposition.
\begin{proposition}\label{prop:rules}
  Given $\pi \in \mathfrak S_n \subset W$, 
  \[
    \pi y_{ij}^k  = y_{\pi(i)\pi(j)}^k.
  \]
  For all $i,j$,
  \[
    g_i y_{ij}^k = \omega y_{ij}^{k-1}, \qquad g_j y_{ij}^k  = y_{ij}^{k+1}.
  \]
  If $\ell \notin \{i,j\}$ then $g_\ell y_{ij}^k = y_{ij}^k$.
\end{proposition}
An elementary computation shows the following result.
\begin{proposition}\label{prop:rules2}
  The product $g_1 g_2 \cdots g_n \in W$ is central in $W$. It acts by
  multiplication by $\omega$ on $K$, and hence multiplication by
  $\omega^{n-1}$ on $Sym^{n-1}(K)$ and on $S_{\G^\perp,-2}$.
\end{proposition}

\section{The internal zonotopal ideal}\label{sec:ideals}
In this section we describe in coordinates the defining ideal of
$S_{\G^\perp,-2}$. We will use this to prove representation stability
in the sense of Church, Ellenberg and Farb and developed for wreath
products with symmetric groups by Gan and Li \cite{ganLi} and Sam and
Snowden \cite{ss}.

\subsection{Generators}\label{sec:generators}
Recall that $\G^\perp \subset K$, where $K$ has coordinates $y_{ij}^k$
with $1 \leq i < j \leq n$ and $0 \leq k \leq m-1$.
\begin{theorem}\label{thm:ideal}
  Let $\G$ be the reflection arrangement of $W = G(m,1,n)$ ($m \geq
  1)$. Then the ideal defining $S_{\G^\perp,-2}$ in $Sym(K)$ is the
  sum of
  \[
    J_1 = \langle y_{ij}^k y_{ij}^{k'} : 1 \leq i < j \leq n, 0 \leq k,k'
    < m \rangle
  \]
  with the smallest $W$-stable ideal $J_2 \subset Sym(K)$ containing
  the single element
  $y_{ij}^0y_{jk}^0 + y_{ik}^0y_{kj}^0 + y_{ji}^0y_{ik}^0$. (Note that $J_1$ is the smallest $W$ stable ideal containing the elements $y_{12}^k y_{12}^{k'}$, $0 \leq k,k' \leq m-1$.)
\end{theorem}
\begin{proof}
  We let $I$ denote the ideal of $S_{\G^\perp,-2}$. We first show that
  $J_1,J_2 \subset I$. Consider the element
  $y_{ij}^k = h_i - \omega^k h_j - h_{ij}^k \in K$. The hyperplanes in
  $\G^\perp$ that do not contain this vector are $h_i$, $h_j$ and
  $h_{ij}^k$, hence $(y_{ij}^k)^2 \in I$. Similarly, the number of
  hyperplanes in $\G^\perp$ that do not contain
  $y_{ij}^k - y_{ij}^{k'} = (\omega^{k'}- \omega^k)h_j - h_{ij}^k +
  h_{ij}^{k'}$ is $3$, hence $(y_{ij}^k - y_{ij}^{k'})^2 \in I$. This
  proves $J_1 \subset I$.  Next we consider
  $y_{ij}^0 + y_{jk}^0 + y_{ki}^0 = -(h_{ij}^0 + h_{jk}^0 +
  h_{ki}^0)$, which is contained in all but $3$ hyperplanes in
  $\G^\perp$, so its square is in $I$. Subtracting off the squared
  variables from $(y_{ij}^0 + y_{jk}^0 + y_{ki}^0)^2$ shows that
  $J_2 \subset I$.

  We claim that $Sym(K)/(J_1+J_2)$ has the same dimension as
  $Sym(K)/I$ as a complex vector space. The non-zero monomials in
  $Sym(K)/J_1$ are in obvious bijection with graphs on vertex set
  $[n]$ whose edges are labeled with the numbers $0,1,\dots,m-1$. We
  will show that the monomials of graphs with cycles are zero in the
  quotient.

  Suppose that we have a monomial in $Sym(K)/(J_1+J_2)$ whose
  corresponding graph has a $3$-cycle. We can use the $W$ invariance
  of this quoient to assume the three cycle has two edges labeled
  $0$. For example, if we had the monomial
  $y_{12}^1 y_{23}^2 y_{31}^3$ then we can apply $g_1^3g_2^{-2}$ and
  work with a multiple of $y_{12}^0 y_{23}^0 y_{31}^6$. Rewriting the
  transformed monomial using the generators of $J_2$ given above, we
  will obtain a sum of two monomials both of which reduce to zero
  modulo $J_1$.

  Suppose now that we have a monomial whose corresponding graph has a
  cycle of length longer than $3$. Similar to the argument above we
  can use the elements in $J_2$ to rewrite this monomial as a sum of
  monomials of graphs with shorter cycles. See
  \cref{fig:cycleShorten}. The resulting sum then reduces to
  zero in the quotient by induction on the length of a cycle.
  \begin{figure}
    \centering
    \includegraphics[scale=.35]{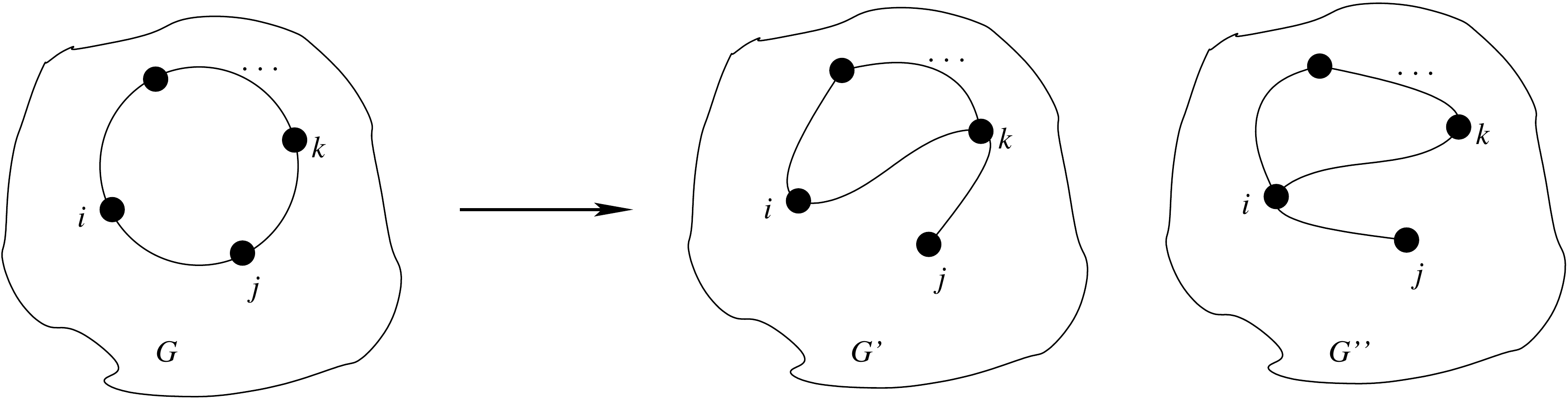}
    \caption{Rewriting a monomial corresponding to a graph with a
      cycle of length larger than $3$ in terms of two monomials whose
      graphs have cycles of smaller length.}
    \label{fig:cycleShorten}
  \end{figure}
  We conclude that $Sym(K)/(J_1+J_2)$ is spanned by monomials in
  bijection with forests on vertex set $[n]$ whose edges are labeled
  with the numbers $0,1,\dots,m-1$.

  Using the relations in $J_2$, we see that out of the forests that
  span $Sym(K)/(J_1+J_2)$, we can use only those where each connected
  component of the underlying forest has its highest label appearing
  as a leaf. This will follow by induction on the degree of the
  highest labeled vertex in each component. Say that the degree of the
  vertex with highest label $v_0$ in a component is at least two, with
  neighbors $v_1$ and $v_2$. Using a relation that swaps out an edge
  $v_0v_1$ or $v_0v_2$ with $v_1v_2$ makes the degree of $v_0$ go down
  by one. Suppose we have a tree that is not a path, where $v_0$ is a
  leaf. Say that $v_1$ is the unique neighbor of $v_0$. If $v_1$ has
  degree larger than two then it has neighbors $v_2$ and $v_3$. The
  relation that swaps an edge $v_1v_2$ or $v_1v_3$ with $v_2v_3$
  decreases the degree of $v_1$.  By induction, a spanning set of
  monomials for the internal zonotopal algebra is indexed by forests
  where each component is a path, and the highest labeled vertex in
  each path is degree one.

  The number such forests with exactly one connected component is at
  once seen to be $m^{n-1}(n-1)!$, which is the dimension of the
  degree $n-1$ piece of $ Sym(K)/I = S_{\G^\perp,-2}$. A calculation
  with the exponential formula shows that the number of such forests
  on $n$ vertices is $1\cdot (m+1)\cdot (2m+1) \cdots ((n-1)m+1)$,
  which is the dimension of $S_{\G^\perp,-2}$.

  We conclude that $\dim Sym(K)/I = \dim Sym(K)/(J_1+J_2)$, and since
  they are both finite dimensional algebras $I = J_1 + J_2$, which is
  what we wanted to show.
\end{proof}

\begin{corollary}\label{cor:cyclic}
  As a representation of $W = G(m,1,n)$, $ S_{\G^\perp,-2}(\textup{top}) $ is
  cyclic with a generator $y_{12}^0y_{23}^0 \dots y_{n(n-1)}^0$. As a
  representation of $G(m,1,n-1)$, $S_{\G^\perp,-2}(\textup{top})$ is
  isomorphic to the regular representation with generator
  $y_{12}^0y_{23}^0 \dots y_{n(n-1)}^0$.
\end{corollary}

\begin{remark}
  There is another way to see that the monomials of graphs containing
  cycles are zero in the quotient. Indeed, if $i_1,i_2,\dots,i_\ell$
  forms a simple cycle then, by definition,
  $(y_{i_1i_2}^0 + y_{i_2 i_3}^0+\dots +y_{i_\ell i_1}^0)^\ell$ is
  zero in the quotient. Expanding, the only square free monomial
  appearing is $y_{i_1i_2}^0 y_{i_2i_3}^0\cdots y_{i_\ell i_1}^0$
  which proves this monomial is in the internal zonotopal ideal. For
  $\ell \leq n-1$ we can act by an appropriate product of the $g_i$s
  to get an arbitrary edge labeled cycle.
\end{remark}

\begin{example}
  Consider the case $n=3$. The monomials in $Sym^2(K)$ that are
  non-zero in $S_{\G^\perp,-2}(\textup{top})$ are the $2^{m+1}$
  monomials of the form
  \[
    y_{12}^ky_{23}^{k'}, \quad
    y_{21}^ky_{13}^{k'},
  \]
  where $0 \leq k,k' < m$ are arbitrary. These monomials form a basis
  for $S_{\G^\perp,-2}(\textup{top})$.
\end{example}

\subsection{Representation Stability}
In this section we write $W_n := G(m,1,n)$. We show that for fixed
degrees $k$ the sequence of representations
$((S_{\G_n,-2})_k)_{n \geq 3}$ of $W_n$ behave in an organized way
explained by the phenomenon of representation stability.

The category $\FI_{\ZZ/m}$ was introduced
by Sam and Snowden in \cite{ss}. Its objects are finite sets and a map
$R \to S$ between two sets is a pair $(f,\rho)$ where $f: R \to S$ is
an injection and $\rho : R \to \ZZ/m$. The composition of
$(f,\rho): R \to S$ with $(g ,\sigma):S \to T$ is defined by
$(g \circ f, \tau)$ where $\tau(x) = \sigma(f(x)) \rho(x)$.

A (complex) representation of $\FI_{\ZZ/m}$ is a sequence
$M= (M_n)_{n \geq 0}$ of finite dimensional complex representations
$M_n$ of $W_n$, together with a sequence of \textit{compatible}
transition maps $M_n \to M_{n+1}$ that are $W_n$-equivariant. Here
\textit{compatible} means that if $w = w_1 \oplus w_2 \in W_{n+\ell}$
where $w_1 \in W_n$, then the action of $w_1$ on $M_n$ agrees with the
action of $w$ on the image of $M_n$ in $M_{n+\ell}$. A simple example
to keep in mind here is the sequence where $M_n$ is the sign
representation of $\mathfrak{S}_n$ for all $n$. The identity maps
$M_n \to M_{n+1}$ do \textit{not} form a compatible sequence (indeed
the sign of $(12)(34)$ is not the same as the sign of $(12)$). On the
other hand, the trivial representation of $\mathfrak{S}_n$ along with
the identity maps \textit{do} form a compatible sequence.

Denote the reflection arrangement associated to the group $W_n$ by
$\G_n$, instead of the usual $\G$.
\begin{proposition}
  For each degree $k \geq 0$, the sequences of representations
  $((S_{\G_n^\perp,-2})_k)_{n \geq 0}$ and
  $((I_{\G_n^\perp,-2})_k)_{n \geq 0}$ are representations of the
  category $\FI_{\ZZ/m}$.
\end{proposition}
\begin{proof}
  Let $K_n$ denote the ambient space of the arrangement
  $\G_n^\perp$. Then $(K_n)_n$ is a $\FI_{\ZZ/m}$ module since given
  $w \in G(m,1,n+\ell)$ which is an extension of an element
  $w_1 \in G(m,1,n)$, we immediately see that the action of $w_1$ on
  $K_n$ agrees with the action of $w$ on the image of $K_n$ in
  $K_{n+\ell}$. 

  It follows that the graded pieces of $Sym(K_n)$ give $\FI_{\ZZ/m}$
  modules. Let $I_n \subset Sym(K_n)$ denote the internal zonotopal
  ideal of $\G_n^\perp$. Restricting the inclusions
  $Sym(K_n) \to Sym(K_{n+1})$ to $I_n$ and apply \cref{thm:ideal} to
  see that the graded pieces of $I_n$ form $\FI_{\ZZ/m}$
  modules. Finally, taking quotients shows that the graded pieces of
  $S_{\G_n^\perp,-2}$ form $\FI_{\ZZ/m}$-modules.
\end{proof}

A representation $M$ of $\FI_{\ZZ/m}$ is said to be finitely generated
if there are finitely many elements in the spaces $M_j$ such that the
smallest $\FI_{\ZZ/m}$ submodule of $M$ containing these elements is
exactly $M$ itself.
\begin{proposition}
  For each degree $k \geq 0$, the sequences of representations
  $((S_{\G_n^\perp,-2})_k)_{n \geq 0}$ and
  $((I_{\G_n^\perp,-2})_k)_{n \geq 0}$ are finitely generated.
\end{proposition}
\begin{proof}
  \cref{thm:ideal} shows finite generation of the graded pieces of the
  sequence of ideals directly. This also follows since the sequence
  $(Sym(K_n)_k)_n$ is finitely generated and the main theorem of
  \cite{ss} proves that subrepresentations of finitely generated
  objects are again finitely generated. The graded pieces of the
  sequence $(S_{\G_n^\perp,-2})_n$ are finitely generated since they
  are quotients of a finitely generated object.
\end{proof}
The importance of being finitely generated here is that it implies
representation stability, as described in
\cite[Definition~1.10]{ganLi}. Recall that the irreducible (complex)
representations of the wreath product $W_n \approx \ZZ/m \wr \mathfrak{S}_n$ are
indexed by $m$-tuples of partitions $\underline \lambda = (\lambda^i)$
satisfying $|\underline \lambda| = \sum_{i=1}^m |\lambda^i| = n$. We
denote this irreducible by $L(\underline\lambda)$; we will not need
the precise construction here. The entries of this $m$-tuple
$\underline \lambda$ are in bijection with the irreducible characters
of $\ZZ/m$ and we assume that $\lambda^1$ corresponds to the trivial
character. For $n$ larger than $|\underline\lambda| + \lambda^1_1$ we
let $\underline \lambda[n]$ denote the partition obtained by adding a
single part of length $n - |\underline \lambda|$ to $\lambda^1$. We
say that a sequence $(M_n)_{n \geq 0}$ of $W_n$ representations is
representation stable if for all $n >\!\! > 0$ the following
conditions hold:
\begin{enumerate}
\item The map $M_n \to M_{n+1}$ is injective.
\item The image of $M_n$ in $M_{n+1}$ generates $M_{n+1}$ as a $W_{n+1}$-module.
\item The irreducible decomposition is given by
  \[
  M_n \approx \bigoplus_{\underline \lambda } L(\underline\lambda[n])^{\oplus m_{\underline\lambda}}
  \]
  for some integers $m_{\underline\lambda} \geq 0$ that do not depend on $n$. 
\end{enumerate}
Gan and Li prove \cite[Theorem~1.12]{ganLi} that $(M_n)_{n \geq 0}$ is
finitely generated if and only if it is representation
stable. Applying their result we have the following.
\begin{theorem}
  For each degree $k \geq 0$, the sequences of representations
  $((S_{\G_n^\perp,-2})_k)_{n \geq 0}$ and
  $((I_{\G_n^\perp,-2})_k)_{n \geq 0}$ are representation stable. 
\end{theorem}

\begin{example}
  One consequence of representation stability of a sequence
  $(M_n)_{n \geq 0}$ is that \textit{$\dim(M_n)$ is eventually
    polynomial in $n$}. In our case the dimension $d_{n,k}$ of
  $(S_{\G_n^\perp,-2})_k$ satisfies the recurrence
  \[
    d_{n,k} = d_{n-1,k} + (n-1)m d_{n-1,k-1}.
  \]
  This follows by computing the coefficient of $q^k$ in the Hilbert
  series,
  \begin{align*}
    (1+mq)(1+2mq) \dots (1+(n-1)mq),
  \end{align*}
  in two ways.  For fixed $k$, $d_{n,k}$ is a polynomial in $n$ for
  every $n \geq 0$.
\end{example}

\section{Type A}\label{sec:typeA}
In this section we investigate the case when $m=1$ in $G(m,1,n)$, so
that our reflection group is the symmetric group of $n$-by-$n$
permutation matrices. To emphasize the dependence on $n$ here we will
write $\mathfrak{S}_n$ instead of $W$.  The main theorem in this
section, \cref{thm:Aind}, will be crucial in our proof of the main
theorem of this paper, \cref{thm:ind}.

\subsection{Statement of the type A result} The reflection arrangement of $\mathfrak{S}_n$ is the well-studied
{braid arrangement} $\A \subset \CC^n$ whose defining hyperplanes are
\[
x_j -x_i,\quad (1 \leq i < j \leq n).
\]
The Macaulay inverse system of the \textit{central} zonotopal ideal $I_{\A,-1}$
has dimension $n^{n-2}$, and was studied by the author and Rhoades in
\cite{bergetRhoades}. There it was shown to be a representation of the
symmetric group $\mathfrak S_n$ that restricted to the well-studied
\textit{parking representation} of $\mathfrak S_{n-1}$. This is the
representation with basis give by sequences
$\mathbf{p}=(p_1,\dots,p_{n-1})$ whose non-decreasing rearrangement
$\mathbf{q}$ satisfies $q_j \leq j$ for all $1 \leq j \leq n-1$.

To identify the Gale dual of $\A$ we label the hyperplane
$x_i - x_j = 0$ by $h_{ij}$. Those $h_{ij}$ with $1 \leq i < j \leq n$
form a basis for $\CC^\A$, and we stipulate that $h_{ji} = -
h_{ij}$. The kernel $K$ of the natural map $\CC^\A \to \CC^n$ is
$\binom{n}{2} -n = \binom{n-1}{2}$ dimensional, and a basis is given
by
\[
  y_{ij} = h_{ij} - h_{in} - h_{jn} \qquad (1 \leq i < j \leq n-1).
\]
One can check that $y_{ji} = -y_{ij}$. This kernel can be naturally
identified with the cycle space of the complete graph on $n$ vertices,
which is spanned by the characteristic vectors of length $3$ oriented
cycles. A basis is then given by those $3$ cycles that visit the
vertex labeled $n$.

The transposition $(kn) \in \mathfrak{S}_n$ acts on the variables
$y_{ij}$ by the rule,
\begin{align}\label{eq:extendedAction}
  (kn) y_{ij} =
  \begin{cases}
    y_{ij} + y_{jk} + y_{ki}, & \textup{if } k \notin \{i,j\},\\
    -y_{ij} & \textup{if } k \in \{i,j\}.
  \end{cases}
\end{align}

Notice that the action of $\mathfrak{S}_{n-1} \subset \mathfrak{S}_n$
on $K$ agrees with action set forth in \cref{sec:galedualG}
with $m=1$ and $n$ replaced by $n-1$.  We have the following analog of
\cref{thm:ind}.
 \begin{theorem}\label{thm:Aind}
   Let $\A$ denote the reflection arrangement of $\mathfrak{S}_n$. Let
   $c \in \mathfrak{S}_{n-1}\subset \mathfrak{S}_n$ denote a any
   $(n-1)$-cycle. As a representation of $\mathfrak{S}_{n-1}$,
   $S_{\A^\perp,-2}(\textup{top})$ is isomorphic to
   $\Ind_{\langle c \rangle}^{\mathfrak{S}_{n-1}}(e^{2\pi i/(n-1)})$.
   As a representation of
   $\mathfrak{S}_{n-2} \subset \mathfrak{S}_{n-1}$,
   $S_{\A^\perp,-2}(\textup{top})$ is isomorphic to the regular
   representation.
 \end{theorem}
 We will prove this result in the next subsection.

 The induced representation here
 $\Ind_{\langle c \rangle}^{\mathfrak{S}_{n-1}}(e^{2\pi i/(n-1)})$ is
 called the Lie representation, which we write as $\Lie_{n-1}$. It is
 the representation of $\mathfrak{S}_{n-1}$ afforded by the
 multilinear component of the free Lie algebra on $n-1$ letters, as is
 shown by the computation of the latter character by Klyachko
 \cite[Corollary~1]{klyachko}. It was shown by Stanley
 \cite[Theorem~7.3]{stanley} that the tensor product of the Lie
 representation with the sign character is isomorphic to the top
 cohomology of the partition lattice, and this is isomorphic to the
 top cohomology of the complement of the braid arrangement
 $\CC^n - \mathcal{A}$.

\subsection{The ideal of $S_{\A^\perp,-2}$ and the Whitehouse representation}
We now compute a generating set for the defining ideal of
$S_{\A^\perp,-2}$. It is possible to do this in the same manner as
$G(m,1,n)$, $m >1$. However, some novel features appear in the
approach we take in this section.

This ideal of $S_{\A^\perp,-2}$ is an ideal in the polynomial ring
$Sym(K) = \CC[y_{ij}: 1 \leq i < j \leq n-1]$.
\begin{lemma}\label{thm:type A ideal}
  The ideal $I_{\A^\perp,-2}$ is equal to
  \[
  J=\la y_{ij}^2 : 1 \leq i < j \leq n \ra + \la y_{ij}y_{ki} +
  y_{ji}y_{kj} + y_{ki}y_{jk} : 1 \leq i < j < k \leq n \ra
  \]  
\end{lemma}
We start by proving one containment.
\begin{proposition}
  The ideal $I=I_{\A_{n+1}^\perp,-2}$ contains $J$.
\end{proposition}
\begin{proof}
  Consider the linear functional $h_{k\ell}^*$ dual to $h_{k\ell}$,
  restricted to $K$. To show that $y_{ij}^2 \in I$ we must compute the
  number of $h_{k\ell}^*$, $k <\ell$, that do not vanish on
  $y_{ij}$. This number is $3$ and hence $y_{ij}^{3-2+1}\in I$.

  The ideal $I$ is stable under the action of $\mathfrak{S}_n$, and
  hence $(kn) y_{ij}^2 = (y_{ij} + y_{jk} + y_{kn})^2\in
  I$. Subtracting off the quadratic terms proves that $I$ contains
  $y_{ij}y_{ki} + y_{ji}y_{kj} + y_{ki}y_{jk}$.
\end{proof}
\begin{proposition}
  The set of polynomials,
  \[
  {y_{ij}^2}, \quad (1 \leq i < j \leq n-1), \qquad
  \underline{-y_{ij}y_{ik}} + y_{ij}y_{jk} - y_{ik}y_{jk},\quad (1
  \leq i < j < k \leq n-1).
  \] 
  forms a Gr\"obner basis under any term order where the leading terms
  are underlined above.
\end{proposition}
\begin{proof}
  We use the fact that syzygies of polynomials with relatively prime
  leading terms need not be computed \cite[Exercise
  15.20]{eisenbud}. By symmetry, the computation reduces to the case
  when $n=4$. This case is easily checked with a computer (e.g., using
  \cite{M2}).
\end{proof}

\begin{proof}[Proof of \cref{thm:type A ideal}]
  The ideal of leading terms of $J$ is
  $\la y_{ij}^2 : 1 \leq i < j \leq n-1\ra + \la y_{ij}y_{jk} : 1 \leq
  i < j < k \leq n-1\ra$. It follows that a basis for $Sym(K)/J$
  consists of square-free monomials in the $y_{ij}$ that avoid
  $y_{ij}y_{ik}$, for $1 \leq i < j < k \leq n$. These monomials are
  in obvious bijection with the set of {decreasing forests} on $n-1$
  vertices. These are forests with vertex set $[n-1]$ where, in every
  component, each path directed away from the largest vertex in that
  component decreases. Adding a vertex labeled $n$ and connecting the
  largest vertex in each component to $n$, we obtain a bijection
  between decreasing forests on $n-1$ vertices, and decreasing trees
  on $n$ vertices. By \cite[Proposition~1.5.5]{EC1}, there are exactly
  $(n-1)!$ such trees.

  It follows that the dimension of $\Sym(K)/J$ as a vector space is
  $(n-1)!$. Since this is the dimension of $S_{\A^\perp,-2}$ (recall,
  its Hilbert series is $(1+q)(1+2q) \cdots (1+(n-2)q)$ and there is a
  containment between $J$ and $I_{\A^\perp,-2}$, the two ideals are
  equal.
\end{proof}
\begin{proof}[Proof of {\cref{thm:Aind}}]
  By \cref{thm:type A ideal}, we have identified
  $S_{\A^\perp,-2}$ with the graded $\mathfrak{S}_{n-1}$ module
  Mathieu studies in \cite[Theorem~3.1]{mathieu}. The top degree
  component of this module is, by his definition, the
  $\mathfrak{S}_{n-1}$ module $\Lie_{n-1}$. By the result of Klyachko
  \cite[Corollary~1]{klyachko} this is the stated induced character.
  Since the induced character is zero at every element not conjugate
  to a power of an $n$ cycle, it is zero at every element of
  $\mathfrak{S}_{n-2}$ except the identity. It follows that the
  restriction of $S_{\A^\perp,-2}(\textup{top})$ to
  $\mathfrak{S}_{n-2}$ is isomorphic to the regular representation.  
\end{proof}

The structure of $S_{\A^\perp,-2}(\textup{top})$ as a representation
of $\mathfrak{S}_n$ is more subtle, but also turns out to have been
studied in work of Mathieu \cite{mathieu} and Gaiffi \cite{gaiffi}. To
state this result we need to define the Whitehouse representation of
the symmetric group $\mathfrak{S}_n$. This is the (\textit{a priori}
virtual) representation,
\[
  Wh_n := \Ind_{\mathfrak{S}_{n-1}}^{\mathfrak{S}_n} (\Lie_{n-1}) - \Lie_n,
\]
which was first studied by Kontsevich in the context of free Lie
algebras and later by Robinson and Whitehouse \cite{whitehouse}.
\begin{theorem}\label{thm:whitehouse}
  There is an isomorphism of representations of $\mathfrak{S}_n$,
  $S_{\A^\perp,-2}(\textrm{top}) \approx Wh_n$.
\end{theorem}

To prove \cref{thm:whitehouse} we will need a result of Sundaram.
\begin{proposition}[Sundaram~{\cite[Lemma~3.1]{sundaram}}]\label{prop:sundaram}
  Let $W_n$ and $V_n$ be (possibly virtual) representations of $\mathfrak{S}_n$,
  and let $\mathbf{1}^\perp$ denote the orthogonal complement of the
  trivial representation in $\CC^n$. Then the following are
  equivalent:
  \begin{enumerate}
  \item $W_n \otimes \mathbf{1}^\perp \approx V_n$.
  \item
    $W_n \oplus V_n \approx
    \Ind_{\mathfrak{S}_{n-1}}^{\mathfrak{S}_n}\Res_{\mathfrak{S}_{n-1}}^{\mathfrak{S}_n}
    W_n$.
  \end{enumerate}
\end{proposition}

\begin{proof}[Proof of \cref{thm:whitehouse}]
  We have identified $S_{\A^\perp,-2}$ with the graded
  $\mathfrak{S}_n$ module Mathieu studies in
  \cite[Theorem~3.3]{mathieu}. We denote the top degree part of this
  $\mathfrak{S}_n$ module by $Q_n$. In \cite[Theorem~4.4]{mathieu} it
  is shown that there is an isomorphism of $\mathfrak{S}_n$-modules,
  \[
    Q_{n} \otimes \mathbf{1}^\perp \approx
    \left(\Res^{\mathfrak{S}_{n+1}}_{\mathfrak{S}_{n}}(Q_{n+1})\right) .
  \]
  Using \cref{thm:Aind} and \cref{prop:sundaram},
  \begin{align*}
    Q_n \oplus \left(\Res^{\mathfrak{S}_{n+1}}_{\mathfrak{S}_{n}}(Q_{n+1})\right) &\approx Q_n \oplus \Lie_n \\
   &\approx \Ind_{\mathfrak{S}_{n-1}}^{\mathfrak{S}_{n}} \left(\Res_{\mathfrak{S_{n-1}}}^{\mathfrak{S_{n}}} Q_n\right) \\
   &\approx \Ind_{\mathfrak{S}_{n-1}}^{\mathfrak{S}_{n}} \left(\Lie_{n-1}\right). 
  \end{align*}
  This proves that
  $Q_n \approx
  \Ind_{\mathfrak{S}_{n-1}}^{\mathfrak{S}_{n}}(\Lie_{n-1}) - \Lie_n$.
\end{proof}

\section{The main theorem} \label{sec:proofs} Our goal in this section
is to prove the main result of the paper.
\begin{theorem}\label{thm:ind}
  Let $\G$ be the reflection arrangement associated to the group
  $W = G(m,1,n)$. Let $C \subset W$ be generated by an $n$ cycle
  $c \in \mathfrak{S}_n \subset W$ and $g_1g_2 \dots g_n \in W$. There is an
  isomorphism of representations,
  \[
    S_{\G^\perp,-2}(\textup{top}) \approx \Ind_{C}^W
    (\chi),
  \]
  where $\chi(c) = e^{2\pi i /n}$ and
  $\chi(g_1 g_2 \dots g_n) = \omega^{n-1}$.
\end{theorem}
Recall that $\omega$ is a primitive $m$th root of unity.
\begin{proof}
  Consider the $\mathfrak{S}_n$ submodule $V_n$ generated by
  $y_{12}^0y_{23}^0 \dots y_{(n-1)n}^0$ in
  $S_{\G^\perp,-2}(\textup{top})$. This module is dimension $(n-1)!$
  by \cref{cor:cyclic} since the monomials of paths where $n$
  is an endpoint form a basis. The quotient,
  \[
S_{\G^\perp,-2}/\langle y_{ij}^k : 1 \leq i \neq j \leq n, k \neq 0 \mod m \rangle
  \]
  is isomorphic as an $\mathfrak{S}_n$ module to
  $S_{\mathcal{A}^\perp,-2}$, where $\A$ is the reflection arrangement
  of $G(1,1,n+1)$ (i.e., the braid arrangement). This follows from the
  description of the ideals in \cref{thm:ideal} and 
  \cref{thm:type A ideal}. The resulting composite
  \[
    V_n  \to S_{\A^\perp,-2}(\textup{top})
  \]
  is an isomorphism of $\mathfrak{S}_n$ modules and hence
  $V_n \approx \Lie_n$.

  By \cref{prop:rules2} and \cref{thm:Aind}, it
  follows that as a representation of
  $\langle\mathfrak{S}_n ,  g_1\cdots g_n \rangle$, 
  $$
  V_n\approx \Ind_{\langle c, g_1 \cdots g_n \rangle}^{\langle\mathfrak{S}_n ,  g_1\cdots g_n \rangle}(\chi),
  $$ 
  where $c$ is any
  $n$-cycle and $\chi$ is the character of defined by
  $\chi (c) = e^{2\pi i/n}$ and $\chi(g_1 \cdots g_n) = \omega^{n-1}$.

  By \cref{cor:cyclic} the $W_n$ orbit of $V_n$ is all of
  $S_{\G^\perp,-2}(\textup{top})$ and so there is a surjective map of
  $W_n$ modules
  \[
    \Ind^{W_n}_{{\langle\mathfrak{S}_n ,  g_1\cdots g_n \rangle}}(V_n) \to S_{\G^\perp,-2}(\textup{top}).
  \]
  The dimension of the module on the left is
  $n!m^n/(n!m) \cdot \dim V_n = (n-1)!m^{n-1}$. Since this is the
  dimension of the module on the right, the map above is an
  isomorphism. By the transitivity of induction \cref{eq:ind}, we have
  \begin{align*}
    \Ind_{\langle c,g_1 \cdots g_n \rangle}^{W_n} (\chi) &\approx
    \Ind^{W_n}_{\langle\mathfrak{S}_n ,  g_1\cdots g_n \rangle} \left(\Ind^{\langle\mathfrak{S}_n ,  g_1\cdots g_n \rangle}_{\langle c,g_1 \cdots g_n \rangle} (\chi)\right)\\
    &\approx \Ind^{W_n}_{\langle\mathfrak{S}_n ,  g_1\cdots g_n \rangle}(V_n)\\ &\approx
    S_{\G^\perp,-2}(\textup{top}).\qedhere
  \end{align*}
\end{proof}

\begin{corollary}\label{cor:coxeter}
If $n$ and $m$ are coprime then 
\[
S_{\G^\perp,-2}(\textup{top}) \approx \Ind_{\langle c \rangle}^W
    (e^{2 \pi i/n} \cdot \omega^{n-1}),
\]
where $c \in W$ is a Coxeter element.
\end{corollary}
A Coxeter element in a well generated reflection group is a product,
in any order, of the generators of the group. Equivalently
\cite{coxeter}, a Coxeter element is an element of $W=G(m,1,n)$ that
has an eigenvalue of multiplicative order $mn$. The multiplicative
order in $W$ of any such element is $mn$.
\begin{proof}
  In the definining representation of $W$, an $n$-cycle $c$ will have
  $e^{2\pi i/n}$ as an eigenvalue and $g_1 \cdots g_n$ has
  $\omega = e^{2\pi i/m}$ as its eigenvalues. Because $(g_1 \cdots g_n)$ is a
  scalar matrix, it follows that $c \cdot (g_1 \cdots g_n)$ has
  $e^{2 \pi i (m+n)/mn}$ as an eigenvalue. 
  Since this is a primative $mn$th root of unity, this element is a
  Coxeter element. Since Coxeter elements in $W$ have order $mn$, it
  follows that the cyclic group appearing in \cref{thm:ind} is
  generated by $c \cdot (g_1 \dots g_n)$ and the character $\chi$
  evaluated at $c \cdot (g_1 \dots g_n)$ is
  $e^{2 \pi i/n} \cdot \omega^{n-1}$.
\end{proof}

\subsection{Factorization of the regular representation} 
One might hope, based on what happens in type A, that
$S_{\G_n^\perp,-2}$ could be made to carry an action of
$G(m,1,n+1)$. However, there does not appear to be an action
generalizing \cref{eq:extendedAction}. The goal of this section is to
prove an analog of the result of Mathieu \cite[Theorem~4.4]{mathieu} occuring in the proof of \cref{thm:whitehouse}, namely,
\[
\Res^{\mathfrak S_{n+1}}_{\mathfrak{S}_n}( Q_{n+1}) 
\approx
Q_n \otimes \mathbf{1}^\perp.
\]
\begin{proposition}\label{thm:factorization}
  Let $E = E_0 \oplus E_1$ be the graded reprentation of $G(m,1,n)$
  that is trivial in degree $0$ and equal to
  $\Ind_{G(m,1,n-1)}^{G(m,1,n)} 1$ in degree $1$. Let $\G_n$ be the
  reflection arrangement of $G(m,1,n)$. Then there is an isomophism of
  graded $G(m,1,n)$-modules,
  \[
   \Res^{G(m,1,n+1)}_{G(m,1,n)}(S_{\G_{n+1}^\perp,-2}) \approx
    S_{\G_n^\perp,-2} \otimes E.
  \]
\end{proposition}
Taking the top degree piece of both sides and applying
\cref{cor:cyclic} we obtain the following result.
\begin{corollary} Maintaining the notation of
  \cref{thm:factorization}, there is an isomorphism of
  $G(m,1,n)$-modules,
\[
\CC[G(m,1,n)] \approx S_{\G_n^\perp,-2}(\textup{top}) \otimes E_1. 
\]  
\end{corollary}
We can unravel this using the proof of Sundaram's
\cref{prop:sundaram}, which brought us to the Whitehouse
representation in \cref{thm:whitehouse}. We obtain a much
more underwhelming result, since we arrive at the simple statement
\[
  \Ind_{G(m,1,n)}^{G(m,1,n+1)} \CC[G(m,1,n) ] \approx \CC[G(m,1,n+1)].
\]
\begin{proof}[Proof of \cref{thm:factorization}]
  From the proof of \cref{thm:ideal}, a linear basis for
  $S_{\G_{n+1}^\perp,-2}$ consists of monomials indexed by certain
  forests on vertex set $[n+1]$. Specifically, the edges are labeled
  with $0,1,\dots,m-1$, each components is a path, and the largest
  labeled vertex in each path has degree one. This means that there is
  a basis indexed by monomials that are either (a) not divisble by any
  variable of the form $y_{(n+1)i}^k$ or (b) divible by exactly one
  variable of the form $y_{(n+1)i}^k$. 

  Let $F = F_0 \oplus F_1$ be trivial in degree $0$ and in degree $1$
  have basis $y_{(n+1)i}^k$, $1 \leq i \leq n$, $0 \leq k \leq
  m-1$. The multiplication map
  $S_{\G_n^\perp,-2} \otimes F \to S_{\G_{n+1}^\perp,-2}$ is
  surjective and is $G(m,1,n)$-equivariant.  Since the Hilbert series
  of both sides match, we get an isomorphism.

  The group $G(m,1,n)$ acts by permuting the variables $y_{(n+1)i}^k$.
  The action is transitive on the set of variables. The stabilizer of
  $y_{(n+1)n}^0$ is precisely $G(m,1,n-1)$ and this proves that
  $E_1 \approx F_1$ as representations of $G(m,1,n)$.
\end{proof}

\section{Type B decreasing trees}\label{sec:typeB}
In this section we compute a G\"robner basis of the ideal
$I_{\G^\perp,-2}$ when $m=2$, where $W$ is the hyperoctohedral group.
We use this to give a non-trivial generalization of the notion of
decreasing trees, as discussed in \cref{sec:typeA}.

We start with a definition. A {$\pm$tree} on $n$ vertices is a
rooted tree with vertex set $[n]$ and root vertex $n$, together with a
$\{+,-\}$-coloring of its edges. A {decreasing $\pm$tree} on
$n$ vertices is a $\pm$tree on $n$ vertices such that on any path
directed away from the root,
\begin{enumerate}
\item along every edge labeled $-$ the vertex labels decrease,
\item  
  for $i<j$ and arbitrary $k$, there is no path of the form 
  \[
    \xymatrix{ j\ar[r]^{-} & i\ar[r]^{+} & k},
  \]
\item for all $i_1< i_2 < i_3 < i_4$, there is no subpath of any of
  the three forms, or their reverse,
  \begin{align*}
    \xymatrix{ i_4 \ar[r]^+ & i_1\ar[r]^{+} & i_2 \ar[r]^+ &i_3},\\
    \xymatrix{ i_4 \ar[r]^+ & i_1\ar[r]^{+} & i_3 \ar[r]^+ &i_2},\\
    \xymatrix{ i_4 \ar[r]^+ & i_2\ar[r]^{+} & i_1 \ar[r]^+ &i_3}.
  \end{align*}
\end{enumerate}
If we ignore the possibility of edges labeled $+$ this definition
reduces to the usual definition of decreasing trees.

\begin{example}
  Here are the $8$ decreasing $\pm$trees on $3$ vertices.
  \[
  \vcenter{\xymatrix@C1px{
     3 \ar[d]^{-}& & 3\ar[d]^+ &&  3\ar[d]^{+} && 3\ar[d]^+ \\
     2 \ar[d]^{-}&&  2\ar[d]^{-} &&  2\ar[d]^{+} && 1\ar[d]^+ \\
     1 &\quad&  1 &\quad& 1  &\quad& 2 
  }}
  \qquad
  \vcenter{
    \xymatrix@C1px{
     &3 \ar[dr]^{-}\ar[dl]_{-} & \\
    1&&2 
  }}
  \quad
  \vcenter{ \xymatrix@C1px{
     &3 \ar[dr]^{+}\ar[dl]_{-} & \\
    1&&2 
  }}\quad\vcenter{
    \xymatrix@C1px{
     &3 \ar[dr]^{-}\ar[dl]_{+} & \\
    1&&2 
  }}\quad\vcenter{
    \xymatrix@C1px{
     &3 \ar[dr]^{+}\ar[dl]_{+} & \\
    1&&2 
  }}
  \]
\end{example}
\begin{example}
  There are $48$ decreasing $\pm$trees on $4$ vertices. Exactly $8$
  such trees are isomorphic to a path rooted at an endpoint, and
  these are displayed below.
  \[
  \vcenter{\xymatrix@C3px{
      4 \ar[d]^{-}& & 4\ar[d]^{+} &&  4\ar[d]^{+} && 4\ar[d]^+ && 4 \ar[d]^{+} && 4\ar[d]^+   &&  4\ar[d]^{+} && 4\ar[d]^+ \\
      3 \ar[d]^{-}& & 3\ar[d]^{-} &&  3\ar[d]^{+} && 3\ar[d]^+ && 3 \ar[d]^{+} && 2\ar[d]^+   &&  2\ar[d]^{+} && 1\ar[d]^{+} \\
      2 \ar[d]^{-}&&  2\ar[d]^{-} &&  2\ar[d]^{-} && 2\ar[d]^+ && 1 \ar[d]^{+} && 3\ar[d]^{+} &&  3\ar[d]^{-} && 3\ar[d]^{-} \\
      1 &\quad& 1 &\quad& 1 &\quad& 1   &\quad&                   2 &\quad&       1 &\quad&       1 &\quad&      2 }}
  \]
\end{example}
\begin{theorem}
  There are $2^{n-1}(n-1)!$  decreasing $\pm$trees on $n$ vertices.
\end{theorem}
\begin{proof}
  We know from the proof of \cref{thm:ideal} that the graphs
  corresponding to monomials not in $I_{\G^\perp,-2}$ are forests on
  $[n]$ whose edges are labeled $0$ and $1$. We change every edge
  labeled $0$ to an edge labeled $-$ likewise with $1$ and $+$. We now
  show that the monomials of decreasing $\pm$trees form a basis for
  the quotient $S_{\G^\perp,-2}(\textup{top})$.

  For this, we claim that a Gr\"obner basis of $I_{\G^\perp,-2}$ is
  furnished by the following monomials:
  \begin{align*}
    (y_{ij}^0)^2,\quad (y_{ij}^0)^2,\quad y_{ij}^0y_{ij}^1,\quad
    y_{ij}^1 y_{jk}^1 y_{ki}^1.
  \end{align*}
  where $i$, $j$ and $k$ range of distinct triples of integers, as well as the polynomials
  \begin{align*}
    \underline{y_{ji}^0y_{ik}^0} + y_{ij}^0y_{jk}^0 + y_{ik}^0y_{kj}^0, \quad    \underline{y_{ji}^1y_{ik}^0} + y_{ij}^1y_{jk}^1 - y_{ik}^0y_{kj}^1,
  \end{align*}
  where the indices range over distinct tuples of integers, and
  finally given $1 \leq i_1 < i_2 < i_3 < i_4 \leq n$ the polynomials,
  \begin{align*}
y^1_{i_1i_2} y^1_{i_1i_3} y^1_{i_2i_4}-y^1_{i_1i_2} y^1_{i_1i_3} y^1_{i_3i_4}-y^1_{i_1i_2} y^1_{i_2i_4} y^1_{i_3i_4}+\underline{y^1_{i_1i_3} y^1_{i_2i_4} y^1_{i_3i_4}},\\
y^1_{i_1i_3} y^1_{i_2i_3} y^1_{i_1i_4}-y^1_{i_1i_3} y^1_{i_2i_3} y^1_{i_2i_4}-y^1_{i_1i_3} y^1_{i_1i_4} y^1_{i_2i_4}+\underline{y^1_{i_2i_3} y^1_{i_1i_4} y^1_{i_2i_4}},\\
y^1_{i_1i_2} y^1_{i_2i_3} y^1_{i_1i_4}-y^1_{i_1i_2} y^1_{i_2i_3} y^1_{i_3i_4}-y^1_{i_1i_2} y^1_{i_1i_4} y^1_{i_3i_4}+\underline{y^1_{i_2i_3} y^1_{i_1i_4} y^1_{i_3i_4}}.
  \end{align*}
  We use any term order where the underlined monomials are leading
  terms (such as graded-reverse lexicographic order with the $y^0$'s
  coming before the $y^1$'s). 

  To see this one applies the Buchberger algorithm to the generators
  of $I_{\G^\perp,-2}$ given in \cref{thm:ideal}. By
  \cite[Exercise 15.20]{eisenbud} we can reduce the computation to the
  case $n=4$, since we need not resolve the syzygies of pair of
  polynomials with relatively prime leading terms. This case is, once
  again, handled at once by a computer to produce the described
  Gr\"obner basis.
\end{proof}

A bijective proof of this result generalizing the one for ordinary
decreasing trees is not known.

\section*{Acknowledgments} I would like to thank an anonymous
referee for their thoughtful review.

\end{document}